\documentclass[12pt]{article}
\usepackage{amsmath}
\usepackage{mathrsfs}
\usepackage{amssymb}
\usepackage{amsfonts,amssymb,amsthm}
\usepackage{color}
\usepackage[all,pdf]{xy}

\allowdisplaybreaks

\newcommand{\Z}{{\mathbb Z}}
\newcommand{\zp}{\Z_+}
\newcommand{\C}{{\mathbb C}}
\newcommand{\G}{\mathcal G}
\newcommand{\U}{\mathcal U}
\renewcommand{\H}{\mathcal H}
\newcommand{\A}{\mathcal A}
\newcommand{\dt}{\delta}
\newcommand{\lmd}{\lambda}
\newcommand{\vf}{\varphi}
\newcommand{\e}{{\epsilon}}
\newcommand{\sgm}{\sigma}

\newcommand{\ing}{\in G}

\newcommand{\hv}{Heisenberg-Virasoro }
\newcommand{\g}{\mathfrak g}
\newcommand{\dd}{\mathfrak d}
\newcommand{\gep}{\g(\e\Z,\lmd)}
\newcommand{\gz}{\g(\Z,\lmd)}
\newcommand{\memu}{M_\e(\mu)}
\newcommand{\mz}{M_\Z(\vf)}

\newcommand{\coef}[2]{{#1}_{(#2)}}
\newcommand{\cli}[1]{C_{LI}^{(#1)}}
\newcommand{\uea}[1]{\mathcal U(\mathfrak{#1})}

\newcommand{\pf}[1]{\noindent{\bf Proof.}#1\hfill{}$\Box$}
\newcommand{\spanc}[1]{\text{span}_{\C}\{#1\}}
\newcommand{\un}[1]{\underline{#1}}
\newcommand{\lm}[1]{L_{-\un{#1}}}
\newcommand{\im}[1]{I_{-\un{#1}}}
\newcommand{\aone}[1]{A_{\un{#1}}}
\newcommand{\atwo}[2]{A_{\un{#1},\un{#2}}}
\newcommand{\vone}[1]{v_{\un{#1}}}
\newcommand{\vtwo}[2]{v_{\un{#1},\un{#2}}}
\newcommand{\mf}[1]{M_{(#1)}}
\newcommand{\iz}[2]{I_{-{#1}_{#2}}}
\newcommand{\lz}[2]{L_{-{#1}_{#2}}}

\newtheorem{thm}{Theorem}[section]

\newtheorem{lem}[thm]{Lemma}

\begin{document}

\begin{center}
{\Large {\bf Verma modules over deformed generalized Heisenberg-Virasoro algebras}}\\
\vspace{0.5cm}
\end{center}

\begin{center}
{Chengkang Xu\footnote{
The author is supported by the National Natural Science Foundation of China (No. 11801375)}\\
Shangrao Normal University, Shangrao, Jiangxi, China\\
Email: xiaoxiongxu@126.com}
\end{center}

\begin{abstract}
Let $\mathfrak g(G,\lambda)$ denote the deformed generalized Heisenberg-Virasoro algebra related to a complex parameter $\lambda\neq-1$ and an additive subgroup $G$ of $\mathbb C$. For a total order on $G$ that is compatible with addition, a Verma module over $\mathfrak g(G,\lambda)$ is defined. In this paper, we completely determine the irreducibility of these Verma modules.
\end{abstract}

\section{Introduction}

\def\theequation{1.\arabic{equation}}
\setcounter{equation}{0}

In this paper we study the irreducibility of Verma modules over deformed generalized \hv algebras, which are generalizations of the \hv algebra, deformed \hv algebras
and generalized \hv algebras.

It is well known that the Heisenberg-Virasoro algebra, first introduced in \cite{ACKP},
is the universal central extension of the Lie algebra of differential operators on a circle
of order no more than one.
The structure and representation theory of the \hv algebra $HV$ has been well developed.
For example, derivations and automorphism group of $HV$ were computed in \cite{SJ}.
Irreducibility for Verma modules over $HV$ was completely determined in \cite{ACKP,B}.
In \cite{LvZ}, L\"{u} and Zhao classified irreducible Harish-Chandra modules over $HV$.
For more, one may see \cite{JJ,ZTL} and references therein.

The Heisenberg-Virasoro algebra $HV$ is graded by $\Z$.
Replacing $\Z$ by an arbitrary additive subgroup $G$ of $\C$,
one gets the so-called generalized \hv algebra $HV[G]$.
Structure and representations of $HV[G]$ were studied in \cite{LG, LZ,SJS}.
In particular, the irreducible Harish-Chandra modules over $HV[G]$
were classified in \cite{LG}.

Infinitesimal deformation of a Lie algebra is one way to give new Lie algebras.
The infinitesimal deformations of the \hv algebra $HV$ are called deformed \hv algebras,
which were given in \cite{LP}.
Inspired by these algebras, Xu \cite{Xu1} introduced deformed generalized \hv algebras $\g(G,\lmd)$,
where $\lmd\neq-1$ is a deformation parameter and $G$ is an additive subgroup of $\C$
such that $G$ is free of rank $\nu$ if $\lmd=-2$.
We also mention that the algebra $\g(G,1)$ is a high rank generalization of the $W$-algebra $W(2,2)$, which was first introduced in \cite{ZD} and extensively studied by others, for example \cite{CG,JP,JZ}.

With respect to a total order on $G$ compatible with addition, there is a triangular decomposition of $\g(G,\lmd)$, and Verma module can be defined accordingly.
In the present paper, we give a complete description of the irreducibility of the Verma modules over $\g(G,\lmd)$ in both cases where the total order on $G$ are dense and discrete. In particular, the irreducibility of the Verma modules over the deformed \hv algebra is determined and used in the proof for the discrete order case for $\g(G,\lmd)$.

We arrange this paper as follows.
In section 2, we introduce the algebra $\g(G,\lmd)$ and its Verma modules.
We prove the irreducibility criterion for Verma modules over the deformed \hv algebra in section 3, and for Verma modules over $\g(G,\lmd)$ in section 4.

Throughout this paper, the symbols $\Z,\zp,\C$ refer to the set of integers,
positive integers, complex numbers respectively. We denote by $\uea m$ the universal enveloping algebra of a Lie algebra $\mathfrak m$. Moreover, any sum is considered as a finite sum, and for a sequence of elements $x_1,\dots,x_n$, we denote by $x_1\cdots \widehat{x_i}\cdots x_n$ the product of these elements with $x_i$ missing.

\section{Verma modules over the algebra $\g(G,\lmd)$}\label{sec2}

\def\theequation{2.\arabic{equation}}
\setcounter{equation}{0}

In this section we give some basics and notations about the deformed generalized \hv algebra $\g(G,\lmd)$ and its Verma modules.

Let $\lmd\neq -1$, and let $G$ be an additive subgroup of $\C$
such that $G$ is free of rank $\nu\geq 1$ if $\lmd=-2$.
Recall from \cite{Xu1} the {\em deformed generalized \hv algebra} $\g(G,\lmd)$
with one deformation parameter $\lmd\neq -1$.
It is formed as the universal central extension of the semi-direct product $LI$ of the centerless generalized Virasoro algebra $L=\spanc{L_a\mid a\ing}$ and its module of intermediate series $I=\spanc{I_b\mid b\ing}$ with action $[L_a,I_b] =(b-\lmd a)I_{a+b}$. When $\lmd=-1$, the algebra $LI$ has no universal central extension, hence the case was not considered in \cite{Xu1}, and also is not considered in this paper. Explicitly, the Lie algebra $\g(G,\lmd)$ has a spanning set 
$\{L_a,I_a, C_L, C_I,\cli i\mid a\ing, 0\leq i\leq\nu\}$ subjecting to Lie brackets
\begin{equation}\label{eq2.1}
 \begin{aligned}
	&[L_a,I_b] =(b-\lmd a)I_{a+b}+\dt_{a+b,0}\left(C_{LI}^{(0)}(a^2+a)\dt_{\lmd,0}
	+\frac1{12}(a^3-a)C_{LI}^{(1)}\dt_{\lmd,1}
	+\sum\limits_{i=2}^\nu a_{(i)}C_{LI}^{(i)}\dt_{\lmd,-2}\right),\\
	&[L_a,L_b] =(b-a)L_{a+b}+\frac1{12}(a^3-a)C_L\dt_{a+b,0},\ \ \ \ \ \ \ \ 
	[I_a,I_b] =aC_I\dt_{a+b,0}\dt_{\lmd,0},
 \end{aligned}
\end{equation}
where $C_L,C_I,\cli i,0\leq i\leq\nu$, are central elements and $\coef ai, 1\leq i\leq\nu$, are coefficients of $a$ with respect to a fixed $\Z$-basis of $\e_1,\cdots, \e_\nu$ of $G$ for the $\lmd=-2$ case. 

For later clarity, we should mention something about the case $\lmd=-2$.
From calculations in \cite{Xu1} the central extension of $LI$ is determined by $\nu-1$ independent nontrivial 2-cocycles.
The term $\sum\limits_{i=2}^\nu a_{(i)}C_{LI}^{(i)}\dt_{\lmd,-2}$ in equation (\ref{eq2.1}) is given by taking the $\nu-1$ nontrivial 2-cocycles with respect to the $\Z$-basis of $\e_1,\cdots, \e_\nu$ of $G$. In particular, when $\nu=1$, i.e., $G\cong\Z$, any 2-cocycle of $LI$ is trivial and the term $\sum\limits_{i=2}^\nu a_{(i)}C_{LI}^{(i)}\dt_{\lmd,-2}$ simply dispears.

When $\lmd=0$, the algebra $\g(G,0)$ is a generalized \hv algebra and the irreducibility of Verma modules over $\g(G,0)$ was determined in \cite{SJS}.
We shall assume that $\lmd\neq 0$ in this paper.
Then equation (\ref{eq2.1}) simply becomes
\begin{equation}\label{eq2.2}
	\begin{aligned}
		&[L_a,I_b] =(b-\lmd a)I_{a+b}+\dt_{a+b,0}\left(\frac1{12}(a^3-a)C_{LI}^{(1)}\dt_{\lmd,1}
		+\sum\limits_{i=2}^\nu a_{(i)}C_{LI}^{(i)}\dt_{\lmd,-2}\right),\\
		&[L_a,L_b] =(b-a)L_{a+b}+\frac1{12}(a^3-a)C_L\dt_{a+b,0},\ \ \ \ \ \ \ \ 
		[I_a,I_b] =0.
	\end{aligned}
\end{equation}

We will simply denote $\g=\g(G,\lmd)$ if no confusion.
Fix a total order $\succ$ on $G$ which is compatible with addition,
i.e., $b\succ c$ implies $a+b\succ a+c$ for any $a,b,c\ing$.
Set $G_+=\{a\in G\mid a\succ 0\}$.
For $a,b\ing$, we shall also write $b\prec a$ if $a\succ b$,
and write $a\succeq b$ if $a\succ b$ or $a=b$.
With respect to the order $\succ$ the algebra $\g$ has a triangular decomposition
$\g=\g_-\oplus\g_0\oplus\g_+$, where
$$\g_0=\spanc{L_0,I_0,C_L,\cli i\mid 1\leq i\leq \nu},\hspace{0.5cm}
  \g_\pm=\spanc{L_a,I_a\mid \pm a\succ0}.$$
Let $\mu$ be a linear function on $\g_0$.
Denote by $I(\mu)$ the left ideal of $\uea g$ generated by
$$\{L_a,I_a,x-\mu(x)\mid a\succ 0, x\in\g_0\}.$$
Then the {\em Verma module with highest weight} $\mu$ over $\g$ with respect to $\succ$
is defined to be
$$M=M(\mu,\succ)=\uea g/I(\mu).$$

For $r>0,\ a_1\succeq a_2\succeq \cdots \succeq a_r\succ 0$,
we write for convenience that $\un a=(a_1,a_2,\cdots,a_r)$
and call $\un a$ a {\em $G_+$-vector}.
Moreover, write
$$\lm{a}=L_{-a_1}\cdots L_{-a_r},\hspace{0.5cm} \im a=I_{-a_1}\cdots I_{-a_r},$$
and denote by $|\un a|$ the length of $\un a$.
Let $\G$ be the set of all $G_+$-vectors.
We may define a total order $\succ$ on $\G$ as follows
(although we use the same symbol for orders on $G$ and $\G$, there should be no ambiguity).
For $\un a,\un b\in\G$, if $k=|\un a|>|\un b|=l$, set $b_{l+1}=\cdots=b_k=0$.
Then we define\\
\centerline{
 $\un a\succ\un b$ if and only if there exists $1\leq i\leq k$ such that $a_i\succ b_i$
  and $a_j=b_j$ for $j<i$.}
Write $v=1+I(\mu)$.
Then by the PBW theorem, the Verma module $M$ has a basis
$$\lm a\im b v,\ \ \ |\un a|,|\un b|\geq 0.$$
We have a subspace filtration of the Verma module $M$
$$0\subset M_0\subset M_1\subset\cdots M_r\subset \cdots\subset M,$$
where $M_r=\spanc{v,\lm a\im b v\mid|\un a|\leq r,\un b\in\G}$ for any $r\geq0$.
We shall write $M_r=0$ if $r<0$.
It is clear that $I_aM_r\subseteq M_{r-1}$ for any $r\in\Z$ and $a\in G_+$.

For $a\in G_+$ set $B_a=\{b\in G\mid a\succ b\succ 0\}$.
The total order $\succ$ on $G$ is called {\em dense} if $B_a$ is infinite for any $a\in G_+$,
called {\em discrete} if $B_a$ is empty for some $a\in G_+$.
Clearly, if $\succ$ is discrete, there is a unique minimal element in $G_+$,
which we will denote by $\e$ in this paper.
Particularly, $\g$ has a subalgebra $\gep$ generated by $\{L_{k\e},I_{k\e}\mid k\in\Z\}$.
Since $\e\Z\cong\Z$ as additive groups, from the statement above for the case $\lmd=-2$ and equation (\ref{eq2.2}), the algebra $\gep=\spanc{L_{k\e},I_{k\e},C_L,\cli1\mid k\in\Z}$ 
subjects to
\begin{equation}\label{eq2.3}
 \begin{aligned}
  &[L_{m\e},L_{n\e}] =(n-m)\e L_{(m+n)\e}+\frac1{12}\left((m\e)^3-m\e\right)C_L\dt_{m+n,0};\ \ \ \
   [I_{m\e},I_{n\e}] =0;\\
  &[L_{m\e},I_{n\e}] =(n-\lmd m)\e I_{(m+n)\e}
       +\frac1{12}\left((m\e)^3-m\e\right)\cli 1\dt_{m+n,0}\dt_{\lmd,1}.
\end{aligned}
\end{equation}
Moreover, the $\gep$-submodule $M_\e(\mu)=\U(\gep)v$ of $M$ is actually a Verma module over $\gep$ with respect to the triangular decomposition $\gep=\gep_-\oplus\gep_0\oplus\gep_+$, where
$$\gep_\pm=\spanc{L_{k\e},I_{k\e}\mid \pm k>0},\ \ \ \
   \gep_0=\spanc{L_0,I_0,C_L,\cli1}.$$
When take $G=\Z$ and one gets the deformed \hv algebra $\gz$ for $\lmd\neq 0,-1$,
which has a basis $\{L_n,I_n,C_L,\cli 1\mid n\in\Z\}$
satisfying
\begin{equation}\label{eq2.4}
	\begin{aligned}
		&[L_m,L_n] =(n-m)L_{m+n}+\frac1{12}(m^3-m)C_L\dt_{m+n,0},\ \ \ \
		[I_m,I_n] =0,\\
		&[L_m,I_n] =(n-\lmd m)I_{m+n}+\frac1{12}(m^3-m)\cli 1\dt_{\lmd,1}\dt_{m+n,0}.
	\end{aligned}
\end{equation}
At last we mention that there is a Lie algebra isomorphism from $\gz$ onto $\gep$ defined by(for $\lmd\neq 0,-1$)
\begin{equation}\label{eq2.5}
	\begin{aligned}
		L_k &\mapsto \e^{-1}L_{k\e}+\dt_{k,0}\frac{\e^{-1}-\e}{24}C_L,\ \ \ \ C_L \mapsto \e C_L,\\
		I_k &\mapsto \e^{-1}I_{k\e}+\dt_{k,0}\dt_{\lmd,1}\frac{\e^{-1}-\e}{24}\cli1,\ \ \ \
		\cli1\mapsto \e \cli1.
	\end{aligned}
\end{equation}

\section{Verma modules over the algebra $\gz$}  \label{sec3}

\def\theequation{3.\arabic{equation}}
\setcounter{equation}{0}

In this section we give a sufficient and necessary condition for a Verma module
over the deformed \hv algebra $\gz$ to be irreducible, which will be used to prove the irreducibility criterion of the Verma module $M(\mu,\succ)$ over $\g(G,\lmd)$ for discrete order $\succ$ in Section \ref{sec4}.

Recall the algebra $\gz$ and its triangular decomposition from the last section.
According to this decomposition one can define a Verma module over $\gz$ as follows.
Let $\vf$ be a linear function on $\gz_0$ and $I(\vf)$ the left ideal of
the universal enveloping algebra $\U(\gz)$ generated by
$$\{L_n,I_n\mid n>0\}\cup\{x-\vf(x)\mid x\in\gz_0\}.$$
Then the {\em Verma module with highest weight }$\vf$ over $\gz$ is defined as the quotient
$$\mz=\U(\gz)/I(\vf).$$
For simplicity denote $\dd=\gz$ and $v=1+I(\vf)$.
By the PBW theorem the Verma module $\mz=\U(\dd_-)v$ has a basis
\begin{equation}\label{eq3.1}
	\iz n1\cdots\iz ns\lz m1\cdots \lz mrv,
\end{equation}
where $r,s\geq0,\ n_1\geq \cdots\geq n_s>0,\ m_1\geq \cdots\geq m_r>0$.
Moreover, $\mz$ has a $\Z$-grading $\mz=\bigoplus_{n\geq 0}\mz_n$
where $\mz_n$ is spanned by vectors of the form in (\ref{eq3.1}) such that
$n_1+\cdots+n_s+m_1+\cdots+m_r=n.$

The algebra $\dd$ has an anti-involution (an anti-isomorphism of order 2) $\sgm$ such that
$$\sgm(L_n)=L_{-n},\ \ \ \sgm(I_n)=I_{-n},\ \ \ \ \sgm(x)=x \text{ for }x\in\dd_0.$$
Moreover, the universal enveloping algebra $\U(\dd)$ has a decomposition
$$\U(\dd)=\U(\dd_0)\oplus\left(\dd_-\U(\dd)+\U(\dd)\dd_+\right).$$
Let $\pi:\uea d\longrightarrow\U(\dd_0)$ denote the projection onto the first summand.
Then we have a symmetric bilinear form $(\cdot \mid\cdot )$ on $\mz$ defined by
$$(xv\mid yv)v=\pi(\widetilde{\sgm}(x)y)v,$$
where $x,y\in\U(\dd_-)$ and $\widetilde{\sgm}$ is the anti-involution of $\uea d$
extended from $\sgm$ by
$$\widetilde{\sgm}(x_1\cdots x_n)=\sgm(x_n)\cdots\sgm(x_1)\ \ \
\text{ for any } x_1,\cdots, x_n\in\dd.$$
Clearly, we have $(v\mid v)=1$ and
$$(xu\mid w)=(u\mid\widetilde{\sgm}(x)w)\ \ \ \text{ for any } x\in\uea d\text{ and } u,w\in\mz.$$
Notice that $(\mz_m\mid\mz_n)=0$ if $m\neq n$,
and the radical of the bilinear form is the maximal $\dd$-submodule of $\mz$.
Then to determine irreducibility of $\mz$,
it suffices to consider the restriction of the bilinear form on each component $\mz_n$.

Define a total order $\succ$ on the set
$\mathcal Z=\{(m_1,\cdots,m_s)\mid s\in\zp, m_1\geq\cdots\geq m_s\in\zp\}$
in the same fashion as the order on $\G$ in Section \ref{sec2}.
Fix an integer $n\geq0$ and denote by $B_n$ the set of the basis of $\mz_n$
consisting of vectors of the form in (\ref{eq3.1}).
We define a total order $\succ$ on $B_n$ as follows. Write
$$\iz n1\cdots\iz nr\lz m1\cdots \lz msv\succ \iz k1\cdots\iz kp\lz l1\cdots \lz lqv$$
if one of the following conditions stands,
\begin{enumerate}
	\item[(C1)] $\sum n_i<\sum k_i$;
	\item[(C2)] $\sum n_i=\sum k_i$ and $(n_1,\cdots,n_r)\succ(k_1,\cdots,k_p)$;
	\item[(C3)] $\sum n_i=\sum k_i,\ (n_1,\cdots,n_r)=(k_1,\cdots,k_p)$ and
	$(m_1,\cdots,m_s)\prec(l_1,\cdots,l_q)$.
\end{enumerate}
Write $B_n=\{u_1,\cdots,u_d\}$ with $u_i\prec u_j$ if $i<j$, where $d=\dim\mz_n$.
Denote $A_n=(A_{ij})$ the $d\times d$ matrix with $A_{ij}=(u_{d+1-i}\mid u_j)$.
In the following we compute the determinant $\det A_n$ of $A_n$.

\begin{lem}\label{lem3.1}
	If $(n_1,\cdots,n_r)\succ(m_1,\cdots,m_s)\in\mathcal Z$, then
	$$(\lz n1\cdots\lz nrv\mid \iz m1\cdots\iz msv)=(\iz m1\cdots\iz msv\mid \lz n1\cdots\lz nrv)=0.$$
\end{lem}
\pf{
	For any integer $m\geq m_1$, we have
	\begin{equation}\label{eq3.2}
		L_m\iz m1\cdots\iz msv=\left(-m(1+\lmd)\vf(I_0)+\frac1{12}(m^3-m)\vf(\cli1)\dt_{\lmd,1}\right)
		\frac\partial{\partial I_{-m}}(\iz m1\cdots\iz ms)v.
	\end{equation}
	Notice that there exists $1\leq k\leq\min\{r,s\}$ such that $n_k>m_k$ and $m_i=n_i$ for $i<k$.
	We see that $L_{n_r}\cdots L_{n_1}\iz m1\cdots\iz msv=0$,
	and then the lemma follows.
}

\begin{lem}\label{lem3.2}
	The determinant $\det A_n$ is a product of a nonzero integer and some
	$$f(k)=-k(1+\lmd)\vf(I_0)+\frac1{12}(k^3-k)\vf(\cli1)\dt_{\lmd,1},\ \ \ 
	k\in\zp.$$
\end{lem}
\pf{
	Let $1\leq a<b\leq d$, then $u_a\prec u_b$. Write
	$$u_a=\iz n1\cdots\iz nr\lz m1\cdots \lz msv,\ \ \ \ \
	u_b=\iz k1\cdots\iz kp\lz l1\cdots \lz lqv.$$
	Then we have
	$$u_{d+1-a}=\iz m1\cdots \iz ms\lz n1\cdots\lz nrv.$$
	Recall the order $\succ$ on $B_n$.
	If case (C1) stands, i.e., $\sum_{i=1}^r n_i<\sum_{i=1}^p k_j$, then we have $\sum_{i=1}^s m_i>\sum_{i=1}^q l_j$.
	It follows from Lemma \ref{lem3.1} that $I_{m_s}\cdots I_{m_1}\lz l1\cdots\lz lqv=0$.
	Hence
	$$L_{n_r}\cdots L_{n_1}\iz k1\cdots\iz kpI_{m_s}\cdots I_{m_1}\lz l1\cdots\lz lqv=0$$
	and $A_{ab}=(u_{d+1-a}\mid u_b)=0$.

	If $\sum_{i=1}^r n_i=\sum_{i=1}^p k_j$, then $\sum_{i=1}^s m_i>\sum_{i=1}^q l_j$ and we have
	$$A_{ab}=(\lz n1\cdots\lz nrv\mid\iz k1\cdots\iz kpv)
	(\iz m1\cdots \iz msv\mid\lz l1\cdots\lz lqv),$$
	which is zero for both case (C2) and case (C3) by Lemma \ref{lem3.1}.
	This proves that the matrix $A_n$ is upper triangular.

	Moreover, by (\ref{eq3.2}) we have
	\begin{align*}
		A_{aa}&=(u_{d+1-a}\mid u_a)
		=(\iz m1\cdots \iz ms\lz n1\cdots\lz nrv\mid\iz n1\cdots\iz nr\lz m1\cdots \lz msv)\\
		&=(\lz n1\cdots\lz nrv\mid\iz n1\cdots\iz nrv)(\iz m1\cdots \iz msv\mid\lz m1\cdots\lz msv)\\
		&=K_a\prod_{i=1}^r f(n_i)^{p_i}\prod_{j=1}^s f(m_j)^{q_j},
	\end{align*}
	where $K_a$ is some nonzero integer,
	$p_i,q_j$ are the times of $n_i,m_j$ appearing in $(n_1,\cdots,n_r)$, $(m_1,\cdots,m_s)$
	respectively.
	This proves the lemma.
}

We can prove our main theorem in this section.
\begin{thm}\label{thm3.3}
	Let $\lmd\neq 0,-1$. The Verma module $\mz$ over $\gz$ is irreducible if and only if $12(1+\lmd)\vf(I_0)-(k^2-1)\vf(\cli1)\dt_{\lmd,1}\neq0$ for any nonzero integer $k$.
\end{thm}
\pf{
	If $12(1+\lmd)\vf(I_0)-(k^2-1)\vf(\cli1)\dt_{\lmd,1}\neq0$ for any nonzero integer $k$, then $f(k)\neq 0$ for any $k\in\zp$.
	Hence the bilinear form on $\mz$ is non-degenerate by Lemma \ref{lem3.2},
	which implies that the $\gz$-module $\mz$ is irreducible.

	Suppose otherwise
	and let $p\in\zp$ be the smallest integer such that
	$12(1+\lmd)\vf(I_0)=(p^2-1)\vf(\cli1)\dt_{\lmd,1}$.
	So $f(p)=0$ and $f(k)\neq 0$ for any $0<k<p$.
	Hence the bilinear form on $\mz$ is degenerate by Lemma \ref{lem3.2},
	whose radical forms a nonzero proper $\gep$-submodule of $\mz$.
}

We mark that for the case $\lmd=1$, $\gz$ is the $W$-algebra $W(2,2)$,
for which the structure of the Verma module $\mz$ is determined in \cite{JZ}.

\section{Irreducibility of Verma modules over $\g(G,\lmd)$}  \label{sec4}

\def\theequation{4.\arabic{equation}}
\setcounter{equation}{0}

In this section we give an irreducibility criterion for the Verma module $M(\mu,\succ)$ over the algebra $\g(G,\lmd)$ with $\lmd\neq 0,-1$, and $G$ being an additive subgroup of $\C$ such that $G$ is free of rank $\nu$ if $\lmd=-2$.
In specific, we have the following theorem.

\begin{thm}\label{thm4.1}
 (1) Suppose the order $\succ$ on $G$ is dense.
     The Verma module $M(\mu,\succ)$ over $\g(G,\lmd)$ is irreducible if and only if
     $\mu(\mathfrak i)\neq 0$, where
     $$\mathfrak i=\spanc{I_0,\cli i\mid 1\leq i\leq \nu}.$$
     When $\mu(\mathfrak i)=0$,
     the Verma module $M(\mu,\succ)$ contains a proper $\g(G,\lmd)$-submodule
     $$N=\spanc{\im bv\mid \un b\in\G},$$
     which is maximal if and only if $(\mu(L_0),\mu(C_L))\neq (0,0)$.
     If $\mu=0$ then $M(\mu,\succ)$ contains a unique maximal $\g(G,\lmd)$-submodule
     $\spanc{\lm a\im bv\mid a,b\in\G,|\un a|+|\un b|>0}$.\\
 (2) Suppose the order $\succ$ on $G$ is discrete.
     The Verma module $M(\mu,\succ)$ over $\g(G,\lmd)$ is irreducible if and only if
     $$24(1+\lmd)\mu(I_0)+\left(1+\lmd+\e^2(1-\lmd-2k^2)\right)\mu(\cli1)\dt_{\lmd,1}\neq0\text{ for any } k\in\Z\setminus\{0\}.$$
\end{thm}
\pf{
(1) Let $M'$ be a nonzero $\g$-submodule of $M$
and recall the subspace $M_0$ of $M$.\\
{\bf Claim 1}: $M'\cap M_0\neq0$. Let
$$u=\sum_{\un a,\un b\in\G}\atwo ab \lm a\im bv\in M',\ \ \ \ \ \atwo ab\neq 0.$$
Notice that since the $L_0$-action on $M$ is semisimple,
we may demand that $u$ is a $L_0$- eigenvector.
This means for all $\un a,\un b$ such that $\atwo ab\neq 0$ in the above equation,
the sum of their entries $a_1+\cdots+a_r+b_1+\cdots+b_s$ are equivalent.

Let $\A=\{\un a\mid \atwo ab\neq 0\text{ for some }\un b\in\G\}$ and
$r=\max\{|\un a|\mid \un a\in\A\}$.
If $r=0$, Claim 1 is trivial.
Assume $r>0$ and denote $\A_r=\{\un a\in\A\mid|\un a|=r\}$.
Since the order $\succ$ is dense, there exists $c\in G_+$ such that $c\prec\min\{a_r\mid \un a\in\A_r\},\ c\notin\{-\lmd a_i\mid \un a\in\A_r, 1\leq i\leq r\}$ and 
$$ \{a_i-c\mid \un a\in\A_r\}\cap\{b_j\mid\atwo ab\neq 0\text{ for all }\un a\in\A_r,1\leq j\leq|\un b|\}=\emptyset.$$
Notice that $I_a M_s\subseteq M_{s-1}$ for any $a\in G_+, s\in \Z$.
We have
\begin{align*}
 I_cu&\equiv \sum_{\un a\in \A_r,\un b\in\G}\atwo ab[I_c,\lm a]\im bv+M_{r-2}\\
     &\equiv-\sum_{\un a\in \A_r,\un b\in\G}\atwo ab\sum_{i=1}^r(c+\lmd a_i)
       L_{-a_1}\cdots\widehat{L_{-a_i}}\cdots L_{-a_r}I_{c-a_i}\im bv+M_{r-2}.
\end{align*}
Since for different $\un a,\un b$, the corresponding vectors
$L_{-a_1}\cdots\widehat{L_{-a_i}}\cdots L_{-a_r}I_{c-a_i}\im bv$
are linearly independent, one gets $I_cu\in M_{r-1}\setminus M_{r-2}$.
Claim 1 follows by induction on $r$.

\noindent {\bf Claim 2}: There exists some vector $\im ev\in M'$.
%for some $\un e=(e_1,\cdots,e_r)\in\G.$
By Claim 1 we may assume the vector $u\in M'$ has the form
$$u=\sum_{\un b\in\G}\aone b\im bv,\ \ \ \aone b\neq 0.$$
Set $T_u=\{\un b\mid \aone b\neq 0\}$ and
let $\un a=(a_1,\cdots,a_r)$ be the maximal element in $T_u$.
We may find some $e_1\in G_+$ such that $e_1\prec a_r,\ a_1+\lmd(a_1-e_1)\neq 0$ and
$$\{x\ing\mid a_1-e_1\prec x\prec a_1\}
   \cap\{b_1,b_2\mid \un b=(b_1,\cdots,b_s)\in T_u\}=\emptyset.$$
Notice that $a_1-e_1-b_j\succ 0$ for any $\un b=(b_1,\cdots,b_s)\in T_u$ and $j\in\{1,\cdots,s\}$
except those $b_j=a_1$.
Then since $u$ is a $L_0$-eigenvector, we get
\begin{equation}\label{eq4.1}
 u_1=L_{a_1-e_1}u=\sum_{\un b}\aone b[L_{a_1-e_1},\im b]v
    =-\sum p(\un b)\aone b(a_1+\lmd(a_1-e_1))I_{-b_2}\cdots I_{-b_s}I_{-e_1}v,
\end{equation}
where the second sum takes over those $\un b\in T_u$ such that $b_1=a_1$,
and $p(\un b)$ is the multiple of $a_1$ in $\un b$.
Let $T_{u_1}$ denote the set of $G_+$-vectors $\un c=(c_1,\cdots,c_s)$ such that
$\{c_1,\cdots,c_s\}=\{e_1,b_2,\cdots,b_s\}$
where $\un b$ appears in the last summand in (\ref{eq4.1}).
One can see that $u_1\neq0$ and
$\un a^{(1)}=(a_2,\cdots,a_r,e_1)$ is the maximal element in $T_{u_1}$.
Set $\aone b^{(1)}=-p\aone b(a_1+\lmd(a_1-e_1))\neq 0$.

For $k=2,\cdots, r$, we define recursively and can easily prove by induction that
\begin{enumerate}
\item[(i)] Let $0\prec e_k\prec e_{k-1},\ a_k+\lmd(a_k-e_k)\neq0$ and
$$\{x\ing\mid a_k-e_k\prec x\prec a_k\}\cap\{b_{k+1},b_{k}\mid\un b\in T_{u_{k-1}}\}=\emptyset.$$
Then
$$u_k=L_{a_k-e_k}u_{k-1}=-\sum p(\un b)\aone b^{(k-1)}(a_k+\lmd(a_k-e_k))
      I_{-b_{k+1}}\cdots I_{-b_s}I_{-e_1}\cdots I_{-e_k}v\neq0,$$
where the sum takes over $\un b\in T_{u_{k-1}}$ such that $b_k=a_k$
and $p(\un b)$ is the multiple of $a_k$ in $\un b$.
\item[(ii)] Let $T_{u_k}$ be the set of $G_+$-vectors formed by
$b_{k+1},\cdots,b_s,e_1,\cdots,e_k$, where $b_{k+1},\cdots,b_s$ appear in $u_k$.
Clearly, $T_{u_k}\neq\emptyset,\
  \aone b^{(k)}=-p(\un b)\aone b^{(k-1)}(a_k+\lmd(a_k-e_k))\neq 0$ and
$\un a^{(k)}=(a_{k+1},\cdots,a_r,e_1,\cdots,e_k)$ is the maximal element in $T_{u_k}$.
\end{enumerate}
Take $k=r$ and one gets Claim 2.

Suppose $\mu(\mathfrak i)\neq 0$, i.e., at least one of $\mu(I_0)$ and $\mu(\cli i)$ is nonzero. Define a $\Z$-linear function $f_\mu$ on $G$ by
$$f_\mu(a)=\mu([L_a,I_{-a}])
  =-(1+\lmd)a\mu(I_0)+\frac1{12}(a^3-a)\mu(\cli1)\dt_{\lmd,1}
  +\sum_{i=2}^\nu \coef ai\mu(\cli i)\dt_{\lmd,-2}.$$
Since the image of $f_\mu$ is not a dense set and the total order $\succ$ on $G$ is dense, we may demand that the vector $\im ev\in M'$ as in Claim 2 satisfies that $f_\mu(e_i)\neq 0$ for all $e_1,\cdots, e_r$. Then we have
$$L_{e_r}\dots L_{e_1}\im ev=\prod\limits_{i=1}^rk_if_\mu(e_i)v\neq 0,$$
where $k_i$ are the multiple of $e_i$ in $\un e$.
This proves $v\in M'$. So $M$ is irreducible.

When $\mu(\mathfrak i)=0$, one can easily check that $N$ is a $\g$-submodule of $M$.
Notice that the quotient $\g$-module $M/N$ is equivalent to the Verma module over
the generalized Virasoro algebra $\spanc{L_a,C_L\mid a\ing}$ generated by $v$.
Then the rest part of (1) follows from Theorem 3.1 (1) in \cite{HWZ}.

(2) Suppose the order $\succ$ on $G$ is discrete and $\lmd\neq 0,-1$.
Recall the minimal element $\e$ in $G_+$, the subalgebra $\gep$ of $\g$
and the $\gep$-module $M_\e(\mu)$.
We first prove that $\g$-module $M(\mu,\succ)$ is irreducible if and only if the $\gep$-module $M_\e(\mu)$ is irreducible.

Write $a\succ\e\Z$ if $a\succ n\e$ for all $n\in\Z$.
Denote $H_+=\{a\ing\mid a\succ\e\Z\},\ H_-=-H_+$ and
$$\g_{H_+}=\spanc{L_a,I_a\mid a\in H_+},\ \ \ \g_{H_-}=\spanc{L_a,I_a\mid a\in H_-}$$
Then we have $G=H_-\cup\e\Z\cup H_+$, $\g_{H_+}M_\e(\mu)=0$ and
$$M=\uea g\otimes_{\U(\gep\oplus\g_{H_+})}\memu=\U(\g_{H_-})\memu.$$
It is clear that the irreducibility of the $\g$-module $M$ implies
the irreducibility of the $\gep$-module $\memu$.

Suppose otherwise that the $\gep$-module $\memu$ is irreducible.
We want to prove that $\memu\cap M''\neq 0$ for any nonzero $\g$-submodule $M''$ of $M$, from which one deduces that the $\g$-module $M$ is irreducible.

Denote by $\H$ the subset of $\G$ consisting of $\un a$ with all entries $a_i\in H_+$.
For $r\geq0$ set
$$\mf r=\spanc{\lm a\im b\memu\mid \un a,\un b\in\H,|\un a|\leq r}+\memu,$$
and for $r<0$ set $\mf r=0$.
Clearly, $I_a\mf r\subseteq\mf{r-1}$ for any $r\in\Z,\ a\in H_+$,
and for any $w\in\mf r$, there exists some $k\in\zp$ such that
$I_{k\e}w\in\mf{r-1}$.\\
{\bf Claim 3}: $\mf0\cap M''\neq0$.
Choose a nonzero vector $u\in M''$.
If $u\in\mf 0$, the claim is trivial.
If $u\in\mf r\setminus\mf{r-1}$ for some $r\in\zp$, then we may write
\begin{equation}\label{eq4.2}
 u=\sum\atwo ab\lm a\im b\vtwo ab+w,
\end{equation}
where the sum takes over $\un a,\un b\in\H$ with $|\un a|=r$ and
$\atwo ab\neq 0,\ 0\neq\vtwo ab\in\memu,\ w\in\mf{r-1}$.
Choose $k\in\zp$ such that $I_{k\e}w\in\mf{r-2},\ I_{k\e}\vtwo ab=0$
for all $\vtwo ab$ in (\ref{eq4.2})
and $\lmd a_i+k\e\neq0$ for all entries $a_i$ of $\un a$ in (\ref{eq4.2}).
Then
$$I_{k\e}u=-\sum\atwo ab\sum_{i=1}^r(\lmd a_i+k\e)L_{-a_1}\cdots
   \widehat{L_{-a_i}}\cdots L_{-a_r}\im bI_{-a_i+k\e}\vtwo ab  \mod\mf{r-2}$$
is a nonzero vector in $M''\cap\mf{r-1}$.
The claim follows by induction on $r$.

By Claim 3 we have a nonzero vector
\begin{equation}\label{eq4.3}
 w=\sum\aone a\im a \vone a\in M'',
\end{equation}
where $\un a\in\H, \aone a\neq 0, 0\neq \vone a\in\memu$.
For $r\geq 0$, denote by $\mf{0,r}$ the subspace of $\mf 0$ spanned by $\im a\memu$
with $\un a\in\H$ and $|\un a|\leq r$,
and for $r<0$ let $\mf{0,r}=0$.
Notice that $\mf {0,0}=\memu$.

Denote by $S_w$ the set of $\un a$ such that $\aone a\neq0$ in (\ref{eq4.3})
and set $r=\max\{|\un a|\mid \un a\in S_w\}$.
If $r=0$ then $w\in M''\cap\memu$ and the proof is done.
Suppose $r>0$ and denote $c=\max\{a_1\mid (a_1\cdots,a_s)\in S_w\}$.
For $k\in\zp$ we have
$$L_{c-k\e} w=-\sum\aone a\sum_{i=1}^{p(\un a)}(a_i+\lmd c-k\lmd\e)
  I_{-a_1}\cdots\widehat{I_{-a_i}}\cdots I_{-a_s}I_{c-a_i-k\e}\vone a,$$
where the first sum takes over $\un a$ with $a_1=c$, and $p(\un a)$ is the times of $I_{-c}$ appearing in $\im a$.
Take $k$ large enough and we see that $L_{c-k\e} w\neq0$ lies in $\mf{0,s}$ for some $s<r$.
Applying induction on $r$ we prove $M''\cap\memu\neq0$.
Till now we have proved that $\g$-module $M(\mu,\succ)$ is irreducible if and only if the $\gep$-module $M_\e(\mu)$ is irreducible.

Recall the algebra $\gz$, its Verma module $\mz$ from Section \ref{sec3}. 
Through the isomorphism given in (\ref{eq2.5}) from $\gz$ onto $\gep$, we consider $M_\e(\mu)$ as a $\gz$-module, which is isomorphic to $\mz$ with $\vf$ satisfying
\begin{align*}
	\vf(L_0) &= \e^{-1}\mu(L_0)+\frac{\e^{-1}-\e}{24}\mu(C_L),
	\ \ \ \ \vf(C_L)=\e\mu(C_L),\\
	\vf(I_0) &= \e^{-1}\mu(I_0)+\frac{\e^{-1}-\e}{24}\mu(\cli1)\dt_{\lmd,1},\ \ \ \
	\vf(\cli1)= \e\mu(\cli1).
\end{align*}
Now applying Theorem \ref{thm3.3}, we get (2).
}

\end{document}